\documentclass[12pt,a4paper]{article}
\usepackage[latin2]{inputenc}
\usepackage{amsmath}
\usepackage{amsfonts}
\usepackage{amssymb}
\usepackage{graphicx}
\usepackage[left=2.00cm, right=2.00cm, top=2.50cm, bottom=2.50cm]{geometry}
\usepackage[T1]{fontenc}

\def\Q{{\mathbb Q}}
\def\Z{{\mathbb Z}}

\newtheorem{lemma}{Lemma}
\newtheorem{theorem}[lemma]{Theorem}

\title{
Calculating power integral bases
in some quartic fields\\
 corresponding to
monogenic families of polynomials
}
\author{
Istv\'an Ga\'al\\
{\small University of Debrecen, Mathematical Institute} \\
{\small H--4002 Debrecen Pf.400., Hungary,} \\
{\small e--mail: gaal.istvan@unideb.hu},
}

\begin{document}
\baselineskip=17pt

\maketitle
\thispagestyle{empty}

\renewcommand{\thefootnote}{\arabic{footnote}}
\setcounter{footnote}{0}

\vspace{0.5cm}


\begin{abstract}
J. Harrington and L. Jones \cite{hj} characterized monogenity of four new
parametric families of quartic polynomials with various Galois groups.
A short time later P. Voutier \cite{vv} added a cyclic family.
In this note we intend to describe all generators of power integral bases in 
the number fields generated by a root of the monogenic polynomials.
\end{abstract}


\section{Introduction}

A number field $K$ of degree $n$ with ring of integers $\Z_K$ 
is called monogenic (cf. \cite{book}) if there 
exists $\xi\in \Z_K$ such that $(1,\xi,\ldots,\xi^{n-1})$ is an integral basis, 
called power integral basis. 
We call $\xi$ the generator of this power integral basis.
$\alpha,\beta\in\Z_K$ are called equivalent,
if $\alpha+\beta\in\Z$ or $\alpha-\beta\in\Z_K$. 
Obviously, $\alpha$ generates a power integral basis in $K$
if and only if any $\beta$, equivalent to $\alpha$ does.
As it is known, any algebraic number field admits up to equivalence only
finitely many generators of power integral bases.

An irreducible polynomial $f(x)\in\Z[x]$ is called monogenic, 
if a root $\xi$ of $f(x)$ generates a power integral basis in $K=\Q(\xi)$.
If $f(x)$ is monogenic, then $K$ is monogenic, but the converse is
not true.

For $\alpha\in\Z_K$ we call the module index $I(\alpha)=(\Z_K:\Z[\alpha])$
the index of $\alpha$. $\alpha$ generates a power integral basis in $K$
if and only if $I(\alpha)=1$.

In a recent paper J. Harrington and L. Jones \cite{hj} gave new infinite
parametric families quartic fields with Galois groups $V_4$ (the Klein four group),
$D_4,A_4,S_4$. They characterized the parameter values for which these
polynomials are monogenic. 
P. Voutier \cite{vv} added a cyclic quartic family.

Our purpose is to describe all generators of
power integral bases in the number fields 
generated by the roots of these monogenic polynomials.
In some cases we succeeded to give a complete description, in some other
cases we made calculations for parameter values in the interval $(-1000,1000)$.
In the following sections we shall refer to the polynomials of \cite{hj} as
$X=2,3,4,5$ as they are denoted in \cite{hj}.

\section{Index form equations in quartic fields}

We briefly present  the results of Ga\'al, Peth\H o and Pohst 
\cite{gppsys}, which allows to reduce index form equations in quartic fields to cubic and quartic Thue equations. 
As there are already efficient tools to solve Thue equations,
which are even implemented in algebraic number theory packages like Magma \cite{magma}, this is an important tool.

Let $K= \Q( \xi )$ be a quartic number field and
$f(x)=x^4+a_1x^3+a_2x^2+a_3x+a_4 \in \Z [x]$ the minimal polynomial of $\xi$.
We represent any $\alpha\in\Z_K$ in the form
\begin{equation}
\alpha \; = \; \frac{a_{\alpha}+x\xi+y\xi^2+z\xi^3}{d}
\label{alfa4}
\end{equation}
with $a_{\alpha},x,y,z\in\Z$, and with a common denominator $d\in\Z$.
Consider the solutions of the equation
\begin{equation}
I(\alpha)= m \;\;\;\;\; (\alpha \in \Z_K )
\label{index4}
\end{equation}
for $0<m\in\Z$. We have 

\begin{theorem} (\cite{gppsys}) Let $i_m=d^6m/n$ where $n=I(\xi)$.
The element $\alpha$ of (\ref{alfa4})
is a solution of (\ref{index4})
if and only if there is a solution $(u,v)\in\Z^2$ of the cubic equation
\begin{eqnarray}
F(u,v)&=&u^3-a_2u^2v+(a_1a_3-4a_4)uv^2\nonumber\\
&&+(4a_2a_4-a_3^2-a_1^2a_4)v^3 = \pm i_m \label{thue34}
\end{eqnarray}
such that $(x,y,z)$ satisfies
\begin{eqnarray}
Q_1(x,y,z)&=&x^2 -xya_1 +y^2a_2+xz(a_1^2-2a_2)+yz(a_3-a_1a_2) \nonumber \\
                &&   +z^2(-a_1a_3+a_2^2+a_4) = u \;\; ,
\nonumber \\
Q_2(x,y,z)&=&y^2-xz-a_1yz+z^2a_2 = v  \;\;\; .
\label{MN4}
\end{eqnarray}
\label{lemma24}
\end{theorem}

For a solution $(u,v)$ of (\ref{thue34}) we set
\[
Q_0(x,y,z)=uQ_2(x,y,z)-vQ_1(x,y,z).
\]
If $\alpha$ in (\ref{alfa4}) is a solution of (\ref{index4}), then
\begin{equation}
Q_0(x,y,z)=0.
\label{q0}
\end{equation}
If $(x_0,y_0,z_0)\in\Z^3$ is a non-trivial solution of (\ref{q0}),
with, say, $z_0\ne 0$ (such a solution can be easily found, see 
L. J. Mordell \cite{mordell}), 
then we can parametrize the solutions $x,y,z$
in the form
\begin{eqnarray}
x&=&rx_0+p\nonumber\\
y&=&ry_0+q\label{pq}\\
z&=&rz_0\nonumber
\end{eqnarray}
with rational parameters $r,p,q$. Substituting these $x,y,z$
into (\ref{q0}) we obtain an equation of the form
\[
r(c_1p+c_2q)=c_3p^2+c_4pq+c_5q^2,
\]
with integer coefficients $c_1,...,c_5$. 
Multiply the equations in (\ref{pq}) by $c_1p+c_2q$
and replace $r(c_1p+c_2q)$ by $c_3p^2+c_4pq+c_5q^2$.
Further multiply the equations in (\ref{pq}) by the square of
the common denominator of $p, q$ to obtain all integer relations
(cf. \cite{gppsim}). We divide those by $\gcd(p, q)^2$ and get
\begin{eqnarray}
kx&=&c_{11}p^2+c_{12}pq+c_{13}q^2\nonumber\\
ky&=&c_{21}p^2+c_{22}pq+c_{23}q^2\label{kpqpq}\\
kz&=&c_{31}p^2+c_{32}pq+c_{33}q^2\nonumber
\end{eqnarray}
with integer $c_{ij}$ and integer parameters $p,q$.
Here $k$ is an integer parameter with the property that 
$k$ divides the $\det(C)/D^2$, 
where $C$ is the 3x3 matrix
with entries $c_{ij}$ and $D$ is the gcd of its entries
(cf. \cite{gppsim}).
Finally, substituting the $x,y,z$ in (\ref{kpqpq}) into
(\ref{MN4}) we obtain
\begin{eqnarray}
F_1(p,q)&=&k^2 u\nonumber\\
F_2(p,q)&=&k^2 v\label{F12}
\end{eqnarray}
\begin{theorem} (\cite{gppsim})
At least one of the equations in (\ref{F12}) is a quartic Thue  equation over the original field $K$.
\label{deg4}
\end{theorem}

Equation (\ref{thue34}) is either trivial to solve (when $F$ is reducible), 
or it is a cubic Thue equation.
Using the above parametrization of the variables $x,y,z$ 
we arrive to the equations (\ref{deg4}), at least one of which
is a quartic Thue equation
over the original number field $K$, see \cite{gppsim}.

\section{X=2: $f_t(x)=x^4+4tx^2+1$, \\ Galois group $V_4$}

The number field generated by a root $\xi$ of $f_t(x)=x^4+4tx^2+1$ is
totally complex for $t\ge 0$ and totally real for $t<0$,
By Theorem 1.1 of J. Harrington and L. Jones \cite{hj},
$f_t(x)$ is monogenic if and only if $4t^2-1$ is square-free.
In this case $(1,\xi,\xi^2,\xi^3)$ is an integral basis in
$K=\Q(\xi)$. Therefore $n=I(\xi)=1$, in (\ref{alfa4}) $d=1$, and we are looking for elements if index $m=1$. Equation (\ref{thue34}) gives
\[
F(u,v)=(u-2v)(u+2v)(u-4tv)=\pm 1.
\]
The system of equations
$u-2v=\pm 1,u+2v=\pm 1,u-4tv=\pm 1$ implies $u=\pm 1,v=0$.

We continue by constructing 
\begin{equation}
Q_1(x,y,z)=x^2+4y^2t-8xzt+z^2(16t^2+1),
\label{q1x2}
\end{equation}
\begin{equation}
Q_2(x,y,z)=y^2-xz+4z^2t.
\label{q2x2}
\end{equation}
By $u=\pm 1,v=0$ we obtain
\begin{equation}
Q_0(x,y,z)=Q_2(x,y,z)=y^2-xz+4z^2t=0,
\label{q0x2}
\end{equation}
having the nontrivial solution $x_0=1,y_0=0,z_0=0$.
Substituting 
\begin{equation}
x=rx_0=r,\;\; y=ry_0+p=p,\;\; z=rz_0+q=q
\label{xyz0x2}
\end{equation}
into (\ref{q0x2})
we obtain $p^2-rq+4q^2t=0$, whence $rq=p^2+4tq^2$.   
Multiplying equations (\ref{xyz0x2}) by $q$ and replacing $rq$
by $p^2+4tq^2$ we obtain
\\
\begin{equation}
\begin{array}{cccc}
kx=&p^2&   &+4tq^2\\
ky=&   &+pq&\\
kz=&&&q^2
\end{array}
\label{kxyzx2}
\end{equation}
where $k=\pm 1$ (see the remarks before Theorem \ref{deg4}).
Substituting $x,y,z$ of (\ref{kxyzx2}) into (\ref{F12}) we obtain
\begin{equation}
F_1(p,q)=p^4+4p^2q^2t+q^4=\pm 1
\label{F1x2}
\end{equation}
and the second equation vanishes.
If $(p,q)\in\Z^2$ is a solution of (\ref{F1x2}), 
then so also is $(- p,- q)$, hence we
only list one of them.

The trivial solutions of (\ref{F1x2}) are $(p,q)=(1,0),(0,1)$.
If $t\ge 0$, then obviously there are no more solutions.
If $t<0$, then we found that for $t=-a^2$ we also have $(1,\pm 2a),(2a,\pm 1)$.
For $-1000\le t<0$  this was proved to be all solutions by solving
(\ref{F1x2}) by Magma. This calculation took 5932 sec.
Substituting the possible pairs $(p,q)$ into (\ref{kxyzx2}) we
obtain all generators of power integral bases.

\begin{theorem}
Assume $4t^2-1$ is square-free. Denote by $K$ the number field
generated by a root $\xi$ of $f_t(x)=x^4+4tx^2+1$.
For any $t$ the elements $\xi$ and $4t\xi+\xi^3$ generate power integral bases.\\
For $t>0$ there are no more inequivalent generators of power integral bases.\\
For $t\le 0$ and $t=-a^2$ ($a\in\Z$) we additionally have $\pm 2a\xi^2+\xi^3$ and
$(1-16a^4)\xi\pm 2a\xi^2+4a^2\xi^3$.
For $-1000\le t$ there are no more inequivalent generators of power integral bases.
\label{thx2}
\end{theorem}

\section{X=3: $f_t(x)=x^4 + 24tx^3 + (12t + 4)x^2 + 4x + 1$,\\
Galois group $D_4$}

The number field $K$ generated by a root $\xi$ of 
$f_t(x)=x^4 + 24tx^3 + (12t + 4)x^2 + 4x + 1$
is totally real for $t<0$, and has two real roots for $t\ge 0$.
By Theorem 1.1 of J. Harrington and L. Jones $f_t(x)$ is monogenic if and only if
$36t^2-1$ is square-free. In these cases (similarly to the case X=2) we have
\[
F(u,v)=(u-12vt)(u^2-4uv+v^2(48t-4)=\pm 1
\]
hence $u-12tv=\pm 1,u=e+12tv$ with $e=\pm 1$. We substitute $u$ into the second factor
of $F(u,v)$  and then obtain two possible cases:\\
A. $(144t^2-4)v^2\pm (24t-4)v+2=0$. This is a second degree equation in $v$.
For $t\ne 0$ the discriminant is negative, for t=0 it is not a square.\\
B: $(144t^2-4)v^2\pm (24t-4)v=0$. In this case either $v=0,u=\pm 1$ or
$v=\pm\frac{24t-4}{144t^2-4}$. For $t\ne 0$ this gives no solution for $v$.
For $t=0$ we obtain $v=\pm 1$, but no solution for $u$.\\

For the solution  $(u,v)=(\pm 1, 0)$ similarly as for $X=2$ we obtain
\[
Q_1(x,y,z)=x^2-24xyt+y^2(12t+4)+xz(576t^2-24t-8)+yz(-288t^2-96t+4)
\]
\[
+z^2(144t^2+17)
\]
\[
Q_2(x,y,z)=y^2-xz-24yzt+z^2(12t+4)
\]
and
\[
Q_0(x,y,z)=Q_2(x,y,z)=y^2-xz-24yzt+z^2(12t+4).
\]
This has the nontrivial solution $x_0=1,y_0=0,z_0=0$. Substituting 
$x=rx_0=r,\;\; y=ry_0+p=p,\;\; z=rz_0+q=q$ into $Q_0(x,y,z)=0$
we obtain $rq=p^2-24tpq+q^2(12t+4)$, whence $rq=p^2-24tpq+q^2(12t+4)$.
Multiplying the equations $x=r,y=p,z=q$ by $q$ and performing the
above substituting of $rq$ we obtain
\begin{equation}
\begin{array}{cccc}
kx=&p^2& -24tpq  &+q^2(12t+4)\\
ky=&   &+pq&\\
kz=&&&q^2
\end{array}
\label{kxyzx3}
\end{equation}
As in X=2 we have $k=\pm 1$. 
Substituting the above representations of $x,y,z$ into 
$Q_1(x,y,z)=u$ and  $Q_2(x,y,z)=v$ we have
\[
F_1(p,q)=p^4-72tp^3q
+p^2q^2(1728t^2+12t+4)+pq^3(-13824t^3-576t^2-192tpq^3+4)
\]
\[
+q^4(6912t^3+2304t^2-96q^4+1)
\]
and $F_2(p,q)$ vanishes. By $p=a+24bt,q=b$ we get
\[
F_1(a,b)=a^4+24ta^3b+(12t+4)a^2b^2+4ab^3+b^4=\pm 1.
\]
We solved this quartic Thue equation by Magma for $-1000\le t\le 1000$.
This calculation took 20476  seconds.
Up to sign we got the solutions 
$(a,b)=(1,0),0,1),(1,-2)$. One can check that these are solutions for any $t$. Moreover
if $t<0$ and $1-2t$ is a square, then we also obtained 
$(a,b)=(1,-1+\sqrt{1-2t}),(1,-1-\sqrt{1-2t})$.
Equivalently, if $t=-(n^2+2n)/2$ is an integer, then
$(a,b)=(1,n),(1,-n-2)$ are solutions.
One can also check that for $t$ of this type, 
these are solutions in general, not only for $-1000\le t$.

By $p=a+24bt,q=b$ and (\ref{kxyzx3}) we obtain

\begin{theorem}
Assume $36t^2-1$ is square-free. Denote by $K$ the number field
generated by a root $\xi$ of $f_t(x)=x^4 + 24tx^3 + (12t + 4)x^2 + 4x + 1$.
For any $t$ the elements \\
$\xi,\;\; (12t+4)\xi+24t\xi^2+\xi^3,\;\; 17\xi+(-2+96t)\xi^2+4\xi^3$
generate power integral bases.\\
Moreover, if  $1-2t$ is a square, or equivalently if $t=-(n^2+2n)/2$
is an integer, then we also have\\
$(-n^4-4n^3+1)\xi+(-2n^4-4n^3+n)\xi^2+n^2\xi^3$ and\\
$(-n^4-4n^3+16n+17)\xi+(-2n^4-12n^3-24n^2-17n-2)\xi^2+(n^2+4n+4)\xi^3$.\\
For $-1000\le t\le 1000$ these are all inequivalent generators of power integral bases.
\label{thx3}
\end{theorem}

\section{X=4: $f_t(x)=x^4 + 2x^3 + 2x^2 + 4tx + 36t^2-16t + 2$,\\
Galois group $A_4$}

The number field $K$ generated by a root $\xi$ of 
$f_t(x)=x^4 + 2x^3 + 2x^2 + 4tx + 36t^2-16t + 2$
is totally complex. 
By Theorem 1.1 of J. Harrington and L. Jones $f_t(x)$ is monogenic
if and only if $(4t-1)(108t^2-54t+7)$ is square-free.
In these cases, similarly to the case X=2 we obtain
\[
F(u,v)=u^3-2u^2v+(72t-144t^2-8)uv^2+(128t^2-64t+8)v^3=\pm 1.
\]
The polynomial $F(u,1)$ has three real roots. 
We solved the above equation by Magma for 
$-1000\le t\le 1000$.
We obtained $(u,v)=(1,0)$ for any $t$ and 
additionally $(u,v)=(1,1),(5,-2)$ for $t=0$.
This calculation took 1600 seconds.
We obtain 
\[
Q_1(x,y,z)=x^2-2xy+2y^2+yz(4t-4)+z^2(-24t+6+36t^2)
\]
\[
Q_2(x,y,z)=y^2-xz-2yz+2z^2.
\]
For the solution $(u,v)=(1,0)$ we set
\[
Q_0(x,y,z)=Q_2(x,y,z)=y^2-xz-2yz+2z^2.
\]
The equation $Q_0(x,y,z)=0$ has a nontrivial solution 
$x_0=1,y_0=0,z_0=0$. Substituting 
$x=rx_0=r,\;\; y=ry_0+p=p,\;\; z=rz_0+q=q$ into $Q_0(x,y,z)=0$ we obtain
$rq=p^2-2pq+2q^2$.
Multiplying by $q$ the equations $x=r,y=p,z=q$ and performing the above
replacement of $rq$ we get
\[
\begin{array}{cccc}
kx=&p^2& -2pq  &+2q^2\\
ky=&   &+pq&\\
kz=&&&q^2
\end{array}
\]
Similarly as before we have $k=\pm 1$. Substituting the above representations of
$x,y,z$ into $Q_1(x,y,z)=\pm 1$ and $Q_2(x,y,z)=0$ we obtain
\[
F_1(p,q)=p^4-6p^3q+14p^2q^2+pq^3(4t-16)+q^4(36t^2-24t+10)=\pm 1
\]
and
$F_2(p,q)$ vanishes.
Substituting $p=a+2b,q=b$ into the first equation we get
\[
F_1(a,b)=a^4+2a^3b+2a^2b^2+4tab^3+b^4(36t^2-16t+2)=\pm 1
\]
We solved this equation by Magma for $-1000\le t\le 1000$, it took
10309  seconds. We found the solution
$(a,b)=(1,0)$ which gives $(p,q)=(1,0)$, $(x,y,z)=(1,0,0)$.

Performing similar calculations in the case $t=0$, for $(u,v)=1,1)$
we get $(x,y,z)=(1,1,0)$ and $(u,v)=(5,-2)$ does not give further solutions.

\begin{theorem}
Assume $(4t-1)(108t^2-54t+7)$ is square-free. Denote by $K$ the number field
generated by a root $\xi$ of $f_t(x)=x^4 + 2x^3 + 2x^2 + 4tx + 36t^2-16t + 2$.
For any $t$ the element $\xi$ generates a power integral basis.
For $t=0$ we additionally have $\xi+\xi^2$.
For $-1000\le t\le 1000$ there are no other inequivalent generators of power integral bases.
\label{thx4}
\end{theorem}

\section{X=5: $f_t(x)=x^4-2x^3-2x^2+6x+4t-2$,\\
Galois group $S_4$}

Let $K$ be the number field generated by a root $\xi$ of
$f_t(x)=x^4-2x^3-2x^2+6x+4t-2$. This polynomial has two real roots for
$t\le 1$ and is totally complex for  $t\ge 2$.
According to Theorem 1.1 of J. Harrington and L. Jones $f_t(x)$
is monogenic if and only if $4t + 1$, $4t -7$ and $64t + 13$ are square-free.
In these cases, similarly to the case X=2 we obtain
\[
F(u,v)=u^3+2u^2v+(-4-16t)uv^2+(-48t-12)v^3=\pm 1.
\]
This polynomial $F(u,1)$ has one real root for
$t\le 1$ and is totally real for  $t\ge 2$.
We solved the above equation for $-1000\le t\le 1000$ by Magma.
It took 39125 seconds. We found the solution
$(u,v)=(1,0)$ for any $t$ and for $t=95$ additionally $(u,v)=(77,-2)$.\\
Corresponding to $(u,v)=(1,0)$ we obtain
\[
Q_1(x,y,z)=x^2+2xy-2y^2+8xz+2yz+z^2(14+4t)
\]
\[
Q_2(x,y,z)=y^2-xz+2yz-2z^2.
\]
Similarly as before we construct
\[
Q_0(x,y,z)=Q_2(x,y,z)=y^2-xz+2yz-2z^2.
\]
A nontrivial solution of $Q_0(x,y,z)=0$ is $x_0=1,y_0=0,z_0=0$.
We set
$x=rx_0=r,\;\; y=ry_0+p=p,\;\; z=rz_0+q=q$, substitute it into $Q_0(x,y,z)=0$
and obtain $rq=p^2+2pq-2q^2$. Multiplying the equations 
$x=r,y=p,z=q$ by $q$ and performing the above substitution of $rq$ we obtain
\[
\begin{array}{cccc}
kx=&p^2& -2pq  &+2q^2\\
ky=&   &+pq&\\
kz=&&&q^2
\end{array}
\]
Similarly as before we have $k=\pm 1$. 
Substituting these representations of $x,y,z$ into $Q_1(x,y,z)=\pm 1$ 
and $Q_2(x,y,z)=0$ we obtain 
\[ 
F_1(p,q)=p^4+6p^3q+10p^2q^2+6pq^3+2q^4+4q^4t=\pm 1
\]
and $F_2(p.q)$ vanishes.
Setting $p=a-2b,q=b$ into the first equation we get
\[
F_1(a,b)=a^4-2a^3b-2a^2b^2+6ab^3-2b^4+4b^4t=\pm 1.
\]
We solved this equation by Magma for $-1000\le t\le 1000$, it took 1671 seconds.
We obtained $(a,b)=(1,0)$ for any $t$ and the following 
solutions for special values of $t$:\\
$t=-736, (a,b)=(7,-1)$\\
$t=-620, (a,b)=(23,3)$\\
$t=-414, (a,b)=(7,1)$\\
$t=-198, (a,b)=(5,-1)$\\
$t=-88, (a,b)=(5,1)$\\
$t=-24, (a,b)=(3,-1)$\\
$t=-6, (a,b)=(3,1)$\\
$t=0, (a,b)=(1,1)$\\
$t=2, (a,b)=(1,-1)$\\
We did not find any rule producing these parameter values or solutions. 

In case $t=95$ we additionally have the solution $(u,v)=(77,-2)$
of $F(u,v)=\pm 1$. Performing the similar calculations we do not 
obtain any solutions for $(a,b)$. Solving 8 quartic Thue equations
took a negligible amount of time by Magma.

\begin{theorem}
Assume $4t + 1$, $4t -7$ and $64t + 13$ are square-free. 
Denote by $K$ the number field
generated by a root $\xi$ of $f_t(x)=x^4-2x^3-2x^2+6x+4t-2$.
For any $t$ the element $\xi$ generates a power integral basis.
For the following parameters we have some additional inequivalent
generators of power integral bases:\\
$t=-736,\;\; 101\xi-9\xi^2+\xi^3$\\
$t=-620, \;\; 205\xi+51\xi^2+9\xi^3$\\
$t=-414,\;\; 17\xi+5\xi^2+\xi^3$\\
$t=-198,\;\; 65\xi-7\xi^2+\xi^3$\\
$t=-88,\;\; 5\xi+3\xi^2+\xi^3$\\
$t=-24,\;\; 37\xi-5\xi^2+\xi^3$\\
$t=-6,\;\; \xi+\xi^2+\xi^3$\\
$t=0,\;\; 5\xi-\xi^2+\xi^3$\\
$t=2,\;\; 17\xi-3\xi^2+\xi^3$\\
For $-1000\le t\le 1000$, up to equivalence, there are no further generators of power integral bases.
\label{thx5}
\end{theorem}

\section{A monogenic family of Galois group $C_4$}

Let $K$ be the number field generated by a root $\xi$ of the polynomial
\[
f := x^4-4tx^3-(1536t^6-512t^4+74t^2-4)x^2-(32768t^9-19456t^7+4608t^5-540t^3+24t)x
\]
\begin{equation} 
-196608t^{12}+163840t^{10}-54784t^8+9216t^6-719t^4+4t^2+2.
\label{cc44}
\end{equation}
P. Voutier \cite{vv} proved that if $16t^2-2$ and $64t^4-64t^2+2$ are squarefree, then 
the polynomial is monogenic, that is $\xi$ generates a power integral basis in $K$.
In this case the polynomial is very complicated, therefore we only have partial results.
However it is interesting to show, how we obtain another generator of power integral basis.
We have
\[
F(u,v)=
u^3+(-4+1536t^6-512t^4+74t^2)u^2v+(80t^2+786432t^{12}-8+716t^4-524288t^{10}+141312t^8
\]\[
-18432t^6)uv^2+(32-145096t^6+17264t^4-1136t^2-8388608t^{12}-813056t^{10}+661504t^8
\]\[
+134217728t^{18}-134217728t^{16}+52953088t^{14})v^3.
\]
The equation $F(u,v)=\pm 1$ obviously has the solution $u=\pm 1,v=0$. 
Then we have $Q_0(x,y,z)=Q_2(x,y,z)$ and 
\[
Q_0(x,y,z)=y^2-xz+4tyz+z^2(-1536t^6+512t^4-74t^2+4)=0
\]
has the nontrivial solution $x_0=1,y_0=0,z_0=0$. Substituting $x=rx_0=r,y=ry_0+p=p,z=rz_0+q=q$ 
into $Q_0(x,y,z)=0$ we obtain 
\[
rq=p^2+4tpq+q^2(-1536t^6+512t^4-74t^2+4).
\]
We multiply the above equations for $x,y,z$ by $q$ and replace $rq$ by the above expression.
Then we obtain
\[
\begin{array}{cccc}
kx=&p^2& -4tpq  &+q^2(-1536t^6+512t^4-74t^2+4)\\
ky=&   &+pq&\\
kz=&&&q^2.
\end{array}
\]
As above we have $k=\pm 1$. Substituting these representation into $Q_1(x,y,z)=1$ we
obtain
\[
F_1(p,q)=
p^4+12tp^3q+(4-26t^2+512t^4-1536t^6)p^2q^2+(7168t^7+12t^3+8t-32768t^9-512t^5)pq^3
\]\[
+(-1024t^6-196608t^{12}+32768t^{10}-1536t^8+2+257t^4-28t^2)q^4=\pm 1.
\]
We substituted $p=a+2tb,q=-b$ to make this formula somewhat simpler. We obtained
\[
G(a,b)=
a^4-4ta^3b+(-74t^2+512t^4+4-1536t^6)a^2b^2+(8t-13312t^7+2560t^5+32768t^9-228t^3)ab^3
\]\[
+(-22016t^8+49t^4-28t^2-196608t^{12}+2+2048t^6+98304t^{10})b^4=\pm 1.
\]
Even this is a quite complicated quartic Thue equation in $a,b$. Solving the equation
for a couple of small parameters by Magma we observed, that in addition to 
the trivial solutions $a=\pm 1,b=0$ for each parameter $t$ some other solutions also occur.
We succeeded to parametrize these values of $a,b$ and we proved in a parametric form, that 
$a=64t^4-12t^2+1,b=4t$ is always a solution. Calculating the corresponding 
$p,q$ and then $x,y,z$,  we obtain

\begin{theorem}
Assume $16t^2-2$ and $64t^4-64t^2+2$ are squarefree and denote by $K$ the number field 
generated by a root $\xi$ of the polynomial (\ref{cc44}).
For any $t$ the inequivalent elements $\xi$ and
\[
(-20480t^8+7680t^6-1104t^4+56t^2+1)\xi+(-256t^5-16t^3-4t)\xi^2+16t^2\xi^3
\]
generate power integral basis.
\end{theorem}

\section{Computational remarks}

All auxiliary calculations were made by Maple \cite{maple}
taking negligible amount of time. The resolution of
Thue equations were performed by using Magma \cite{magma},
requiring all together 79113 seconds, or equivalently 
about 22 hours.

\end{document}